\def\draft{y}
\documentclass{amsart}
\usepackage{fullpage,amssymb,epic,eepic,epsfig}

\theoremstyle{plain}

\newtheorem{theorem}{Theorem}
\newtheorem{proposition}{Proposition}[section]
\newtheorem{lemma}[proposition]{Lemma}
\newtheorem{corollary}[proposition]{Corollary}

\newtheorem{observation}{Observation}

\theoremstyle{definition}
\newtheorem{definition}[proposition]{Definition}

\newtheorem{question}{Question}

\theoremstyle{remark}
\newtheorem{example}[proposition]{Example}

\newtheorem{remark}[proposition]{Remark}

\def\printname#1{
	\if\draft y
		\smash{\makebox[0pt]{\hspace{-0.5in}
			\raisebox{8pt}{\tt\tiny #1}}}
	\fi
}

\newcommand{\psdraw}[2]
         {\begin{array}{c} \hspace{-1.3mm}
	\raisebox{-4pt}{\epsfig{figure=draws/#1.eps,width=#2}}
	\hspace{-1.9mm}\end{array}}

\newlength{\standardunitlength}
\catcode`\@=11
\long\def\@makecaption#1#2{%
     \vskip 10pt

\setbox\@tempboxa\hbox{
       \small\sf{\bfcaptionfont #1. }\ignorespaces #2}%
     \ifdim \wd\@tempboxa >\captionwidth {%
         \rightskip=\@captionmargin\leftskip=\@captionmargin
         \unhbox\@tempboxa\par}%
       \else
         \hbox to\hsize{\hfil\box\@tempboxa\hfil}%
     \fi}
\font\bfcaptionfont=cmssbx10 scaled \magstephalf
\newdimen\@captionmargin\@captionmargin=2\parindent
\newdimen\captionwidth\captionwidth=\hsize
\catcode`\@=12

\def\lbl#1{\label{#1}\printname{#1}}

\def\BN{\mathbb N}

\def\AC{\mathcal AC}

\def\C{\mathcal C}

\def\D{\Delta}
\def\K{\mathcal K}

\def\F{\mathcal F}

\def\Err{\mathcal Err}

\def\l{\lambda}

\def\S{\Sigma}

\def\as{algebraically split}

\def\w{\omega}

\def\e{\epsilon}

\def\d{\delta}
\def\b{\beta}

\def\rot{\mathrm{rot}}
\def\longto{\longrightarrow}
\def\Sev{{\mathcal S}}

\def\sign{\mathrm{sign}}
\def\exc{\mathrm{exc}}

\def\mult{\mathrm{mult}}
\def\as{\mathrm{as}}

\def\longto{\longrightarrow}

\def\defect{\mathrm{def}}
\def\deg{\mathrm{deg}}

\begin{document}


\title[Chromatic Polynomial, Colered Jones Function and q-Binomial Counting]{Chromatic Polynomial,
Colored Jones Function and q-Binomial Counting}

\author{Martin Loebl}
\address{Dept.~of Applied Mathematics and\\
Institute of Theoretical Computer Science (ITI)\\
Charles University \\
Malostranske n. 25 \\
118 00 Praha 1 \\
Czech Republic.}
\email{loebl@kam.mff.cuni.cz}

\thanks{
1991 {\em Mathematics Classification.} Primary  05A30  . Secondary 57N10 .
}

\date{
This edition: \today \hspace{0.5cm} First edition: February 18, 2003.}

\begin{abstract}
We define a q-chromatic function on graphs, list some of its properties and 
provide some formulas in the class of general chordal graphs. Then we relate the q-chromatic function
to the colored Jones function of knots. This leads to a curious expression
of the colored Jones function of a knot diagram $\K$ as a 'defected chromatic operator' 
applied to a power series whose coefficients are linear combinations
 of chord diagrams constructed from 'flows' on reduced $\K$. 
\end{abstract}

\maketitle

\tableofcontents


\section{Introduction and Statement of the Results}
\lbl{sec.int}

The main purpose of this paper has been a desire to recast the complicated combinatorial
construction of 'categorification of flows' of \cite{GL}
in more common combinatorial terms. This lead to a definition of the q-chromatic
function which I have observed started to live by itself.
The final formula may be related to the Kontsevich integral, and perhaps also
to the Khovanov's work on categorification of the Jones polynomial (\cite{Kh}).
Final push to finish first version of the manuscript came from Stavros Garoufalidis, who informed me about paper \cite{HGR} which establishes a Khovanov homology theory for the chromatic polynomial of graphs.  

Let me sketch main steps of the treatment of the colored Jones function $J_n$.
We start with the state sum of the $n$-cabling of a knot. Let us denote by $w(s)$
the contribution of state $s$. First we associate to each state $s$ a triple
$(f(s), S(s), v(s))$, where $f(s)$ is a non-negative integer flow; 
you can imagine that the flow lives in the reduced knot diagram $K$, eventhough it turns out
to be more convenient to define it on the {\em arc graph of the knot diagram}.

$S(s)$ is a set system on the set of $\sum_ef(s)(e)$ elements where the sum is over 
all 'jump-up' transitions $e$ of $K$ (which are later associated to 'red edges' of the arc
 graph).  Each such set system will be called simply '$f(s)$-structure'.
The number of $f(s)$-structures for a given flow $f(s)$ is given by a product of binomial
coefficients. 

Finally $v(s)$ is a non-negative integer vector of length $\sum_ef(s)(e)$, and for each $i$, $0\leq v(s)_i\leq (n-1)$.

We then realise that we can naturally represent $S(s)$ as a chord diagram $ChD(s)$
with $\sum_ef(s)(e)$ chords.

If we denote by $w(n,f,S)$ the sum of $w(s)$ over all states $s$ with $f(s)= f$, $S(s)= S$
and by $G(s)$ the intersection graph of the chords of $ChD(s)$ then
we observe that $w(f,S)$ may be written as
 
$$
Z_{n,f}(t)M_t^{\defect}(G(s),n),
$$
where $M_t^{\defect}$ is a 'defected' q-chromatic function and $Z_{n,f}(t)$ is a Laurent
polynomial in $t$ whose precise form is given in Theorem \ref{thm.ma2}.

This leads to a curious expression
of the colored Jones function of a knot diagram $\K$ as a power series whose coefficients are 
equal to a 'defected chromatic operator' applied to linear combinations
 of chord diagrams constructed from 'flows' on reduced $\K$. 

A graph is a pair $G=(V,E)$ where $V$ is a finite set of {\it vertices} and $E$ is 
 a set of unordered pairs of elements of $V$, called {\it edges}. 
 If $e=xy$ is an edge then the vertices $x,y$ are called {\it end-vertices} of $e$. 
A graph $G'=(V',E')$ is called a {\it subgraph} of a graph $G=(V,E)$ if
 $V'\subset V$ and $E'\subset E$.

\subsection{q-Bichromate}
\lbl{sub.qchrr}
In this paper we study the following function on graphs:

\begin{definition}
\lbl{def.qc}
Let $G=(V,E)$ be a graph. Let $V=\{1,\dots, k\}$ and let $V(G,n)$ denote the set of all vectors 
$(v_1,\dots, v_k)$ such that $0\leq v_i\leq n-1$ for each $i\leq k$ and $v_i\neq v_j$
whenever $\{i,j\}$ is an edge of $G$. We let 
$$
M_q(G,n)=\sum_{(v_1\dots v_k)\in V(G,n)}q^{\sum_i v_i}.
$$
\end{definition}

Note that $M_q(G,z)|_{q=1}$ is the classic chromatic polynomial of $G$.
If $G=(V,E)$ is a graph and 
$A\subset E$ then let $C(A)$ denote the set of the connectivity components of
graph $(V,A)$, and if $W\in C(A)$ then let $|W|$ denote the number of vertices
of $W$. 

\begin{theorem}
\lbl{thm.qch}
$$
M_q(G,z)=\sum_{A\subset E} (-1)^{|A|} \prod_{W\in C(A)} (z)_{q^{|W|}}.
$$
 \end{theorem}

The following function extensively studied in combinatorics is called {\it bichromate}:
$$
B(G,a,b)= \sum_{A\subset E} a^{c(A)}b^{|A|}.
$$

Note that the bichromate is equivalent to the Tutte polynomial (see next section
for more details).
 
For $n>0$ let $(n)_q= \frac{q^n-1}{q-1}$ be a quantum integer.  
We let $(n)!_q=\prod_{i=1}^n (i)_q$ and for $0\leq k\leq n$ we define 
the quantum binomial coefficients by
$$
{{n}\choose{k}}_q=\frac{(n)!_q}{(k)!_q(n-k)!_q}.
$$

The formula of Theorem \ref{thm.qch} leads naturally to a definition of 
{\it q-bichromate}.

\begin{definition}
\lbl{def.qt}
We let
$$
B_q(G,x,y)= \sum_{A\subset E} x^{|A|} \prod_{W\in C(A)} (y)_{q^{|W|}}.
$$
\end{definition}

Note that $B_{q=1}(G,x,y)= B(G,x,y)$.

It is well known that the bichromate counts several interesting things
in statistical physics. We concentrate on the Potts and 
Ising partition functions, and on the Jones polynomial, and discuss their 
q-extensions. 

\subsection{Potts partition function}
\lbl{sub.p}

\begin{definition}
\lbl{def.P}
Let $G= (V,E)$ be a graph, $k\geq 1$ integer and $J_e$ a weight (coupling constant)
associated with edge $e\in E$. The Potts model partition function is defined as
$$
P^k(G, J_e)= \sum_{s}e^{E(P^k)(s)},
$$
where the sum is over all functions $s$ from $V$ to $\{1,\dots, k\}$ and 
$$
E(P^k)(s)= \sum_{\{i,j\}\in E} J_{ij}\delta(s(i), s(j)).
$$
\end{definition}

Following \cite{VN}, we may write 
$$
P^k(G, J_e)= \sum_s\prod_{\{i,j\}\in E}(1+ v_{ij}\delta(s(i), s(j)))= 
\sum_{A\subset E} k^{c(A)} \prod_{\{i,j\}\in A}v_{ij},
$$
where $v_{ij}= e^{J_{ij}}-1$. If all $J_{ij}$ are the same, we get an expression
of the Potts partition function in the form of the  bichromatic polynomial:

\begin{theorem}
\lbl{thm.ptts}
$$
P^k(G, x)= \sum_{s}\prod_{\{i,j\}\in E}e^{x\delta(s(i), s(j))}= 
\sum_{A\subset E} k^{c(A)} (e^x-1)^{|A|}= B(G, e^x-1, k).
$$
\end{theorem}

\subsection{q-Potts}
\lbl{sub.qptts}

What happens if we replace $B(G, e^x-1, k)$ by $B_q(G, e^x-1, k)$? It turns out that
this introduces an additional external field to the Potts model.

\begin{theorem}
\lbl{thm.qP}
$$
\sum_{A\subset E} \prod_{W\in C(A)} (k)_{q^{|W|}}
\prod_{\{i,j\}\in A}v_{ij}=
\sum_{s}q^{\sum_{v\in V}s(v)}e^{E(P^k)(s)},
$$
where $v_{ij}= e^{J_{ij}}-1$ as above.
\end{theorem}

\subsection{Ising partition function}
\lbl{sub.i}

The Ising partition function $Z(G)$ of a graph $G$ is equivalent to $P^2(G)$: 
$$
Z(G, J_e)= \sum_{s}e^{E(Z)(s)},
$$
where the sum is over all functions $s$ from $V$ to $\{1,-1\}$ and 
$$
E(Z)(s)= \sum_{\{i,j\}\in E} J_{ij}s(i)s(j).
$$

We immediately have 
$$
Z(G, J_e)= e^{-\sum_{\{i,j\}\in E} J_{ij}}P^2(G, 2J_e).
$$

Not surprisingly, q-bichromate again adds an external field to the Ising partition function.

\begin{corollary}
\lbl{cor.qZ}
$$
Z(G, x)= e^{-|E|x} B(G, e^{2x}-1, 2) 
$$
and
$$
\sum_{s}q^{\sum_{v\in V}s(v)}e^{E(Z)(s)}= q^{-3|V|x} e^{-|E|x} B_q(G, e^{2x}-1, 2) .
$$
\end{corollary}

\subsection{Van der Waerden Theorem}
\lbl{sub.m}

A remarkable fact about the Ising partition function is a theorem of Van der Waerden
which expresses it using the generating function of even subgraphs. It is not hard
to formulate its q-generalisation.
We use the following notation:

$$
sinh(x)=\frac{e^x-e^{-x}}{2}, cosh(x)=\frac{e^x + e^{-x}}{2}, 
 th(x)=\frac{sinh(x)}{cosh(x)}.
$$

\begin{theorem}
\lbl{thm.w}
$$
\sum_{s}q^{\sum_{v\in V}s(v)}e^{\sum_{\{i,j\}\in E} J_{ij}s(i)s(j)}=
$$
$$
\prod_{\{i,j\}\in E}cosh( J_{ij})\sum_{A\subset E}\prod_{\{i,j\}\in A}th( J_{ij})
(q- q^{-1})^{o(A)}(q+ q^{-1})^{|V|- o(A)},
$$
where $o(E)$ denotes 
the number of vertices of $G$ of an odd degree.
\end{theorem}

\subsection{Jones polynomial}
\lbl{sub.q}

Following \cite{VN} we show how the Jones polynomial may be derived from the Potts partition 
function. Let $G= (V,E)$ be a planar directed graph which is
a directed knot diagram; hence each vertex (crossing) $v$ has a sign $\sign(v)$ 
associated with it, and two arcs entering and leaving it.

Given $G$, we construct its {\it median graph} $M(G)= (V(M(G)), E(M(G)))$ as follows: 
color the faces of $G$ white and black so that neighbouring faces receive a different color. 
Assume the outer face is white; then let $V(M(G))$ be the set of the black faces,
and two vertices are joined by an edge if the corresponding faces share a crossing.
Note that $M(G))$ is again a planar graph. For edge $e$ of $M(G)$ let $b(e)$ be the sign of 
the crossing shared by the end-vertices of $e$. 

Let us now describe what will a {\it state} be: we can 'split' each vertex $v$ of $G$
so that the white faces incident with $v$ are joined into one face and the black faces are 
disconnected, or vice versa. Let $|V|= n$.
There are $2^n$ ways to split all the vertices of $G$: the ways are called {\it states}.
After performing all the splittings of a state $s$, we are left with a set of disjoint
non-self-intersecting cycles in the plane; let $S(s)$ denote their number.
For vertex $v\in V$ let $\e_v(s)= 1$ if $s$ splits $v$ so that the black faces are joined
(i.e. the corresponding edge of the median graph is not cut) and let $\e_v(s)= -1$ otherwise. 
 The following statement consists of Theorems 2.6 and 2.8 by Kauffman (see \cite{K1}).

\begin{theorem}
\lbl{thm.K1}
Let $G= (V, E)$ be an oriented knot diagram. The following function $f_K(G)$ is a knot 
invariant:
$$
f_K(G,A)= (-A)^{-3W(G)}\sum_s (-A^2-A^{-2})^{S(s)-1} 
A^{\sum_{v\in V}\sign(v)\e_v(s)},
$$
where $W(G)= \sum_{v\in V} \sign(v)$. Moreover the Jones polynomial equals
$$
J(G, A^{-4})= f_K(G, A).
$$
\end{theorem}

Each state $s$ determines a subset of edges $E(s)$ of $M(G)$, which are not cut
by the splittings of $s$, and it is easy to see that this gives a natural bijection between the 
set of states and the subsets of edges of $M(G)$. Moreover  for each state $s$
$$
\sum_{v\in V}\sign(v)\e_v(s) =  -\sum_{e\in E(M(G))} b(e) + 2\sum_{e\in E(s)}b(e).
$$ 

\begin{proposition}
\lbl{prop.mm}
$$
S(s)= 2c(E(s))+ |E(s)|- |V(M(G))|,
$$
where $c(E(s))$ denotes the number of connectivity components of $(V(M(G)),E(s))$.
\end{proposition}
\begin{proof} 
Note that $S(s)= f(E(s))+ c(E(s)) -1$, where $f(E(s))$ denotes the number of faces
of $(V(M(G)), E(s))$. Hence the formula follows from the Euler formula for the planar graphs.
\end{proof}

\begin{corollary}
\lbl{cor.kk}
$$
f_K(G,A)= (-A)^{-3W(G)}(-A^2-A^{-2})^{-|V(M(G))|-1}
A^{-\sum_{e\in E(M(G))} b(e)}
 \sum_s (-A^2-A^{-2})^{2c(E(s))} \prod_{e\in E(s)}(-A^2-A^{-2})A^{2b(e)}.
$$
\end{corollary}

This provides an expression of the Jones polynomial of an arbitrary knot-diagram $G$ with the same
sign $b$ of each crossing as a bichromate.

\begin{corollary}
\lbl{cor.kkk}
$$
J(G,A^{-4})= (-A)^{-3V(M(G))}(-A^2-A^{-2})^{-|V(M(G))|-1}A^{-b|E(M(G))|}
B(M(G), (-A^2-A^{-2})A^{2b}, (-A^2-A^{-2})^2).
$$ 
\end{corollary}

\begin{question} 
Is $B_q(M(G), (-A^2-A^{-2})A^{2b}, (-A^2-A^{-2})^2)$ times an appropriate constant
 also a knot invariant?
\end{question}

We remark that the Jones polynomial of an alternating link also equals to a specialization of the 
{\em Tutte polynomial} of the median graph of its planar projection, \cite[Proposition 5.2.14]{W}. 

\subsection{Aproximating the Jones polynomial}
\lbl{sub.aprr}

A {\em positive link} (resp. negative link) is one that has a planar 
projection with positive (resp. negative) crossings only. Notice that
the mirror image of a positive link is negative, and vice versa.
As remarked by Stavros Garoufalidis, positive links and positive braids play an
important role in symplectic aspects of smooth 4-dimensional topology
(such as existence of Lefsetz fibrations). For a
discussion of how well-known invariants of links behave when
restricted to the class of positive links, see \cite{CM,St}.

It is well-known that computing the Jones polynomial is a \#P-{\em complete}
problem (see \cite[Sec.6]{W}). On the other hand, one may ask about {\em approximating} the
Jones polynomial. Partial results on the existence of {\em Fully polynomial randomized 
approximation scheme} (FPRAS, in short) that approximate values of the Tutte polynomial 
are known, see \cite{AFW}. Hence, it is still possible that Jones polynomial for alternating
and positive links may be well approximable.

\subsection{q-Chromatic function of chordal graphs}
\lbl{sub.chordal}

Next we study chordal graphs, i.e. graphs such that each cycle of
length at least four has a chord. Let $G=(V,E)$
be a chordal graph. We fix a linear ordering $x_1,\dots, x_k$ of the vertices
of $V$ so that for each $i$, vertex $x_i$ is a simplicial vertex, i.e. its
neighbourhood is complete, in the subgraph induced by vertices $x_1,\dots, x_i$;
let $m(i)$ denote the number of vertices in that neighbourhood of $x_i$. It is well-known
that the existence of such an order of vertices characterises the chordal graphs.

A tree is an acyclic connected graph. We consider trees {\it rooted}, i.e. a vertex $r$ is 
distinguished in each tree. Hence we will denote trees by a triple $T= (V,E,r)$. For each
vertex $x\neq r$ of $T$ there is unique path in $T$ connecteing it to $r$. 
The neighbour of $x$ on it is called {\it predecessor of $x$} and denoted by $p(x)$. The set of 
 the vertices of the path without $x$ is denoted by $P(x)$. We say that a subtree of a rooted
tree {\em starts} at its unique nearest vertex to the root. 

Another well known characterisation says that a graph is chordal if and only if
it is an intersection graph of subtrees of a tree.
Let $G=(V,E)$ be a chordal graph and $x_1,\dots, x_k$ the specified ordering of its vertices. Let $T=(W,F)$ be the tree whose subtrees 'represent' $G$, i.e. there
are subtrees $T_v, v\in V$ of $T$ so that $T_u\cap T_v\neq \emptyset$ if and only if
$uv$ is an edge of $G$. We can choose
a vertex $r\in W$ arbitrarily as a root of $T$ and define sets $A_w, B_w, w\in W$
as follows: $A_w=\{v; T_v \text{starts at}  w\}$ and 
$B_w=\{v; T_v \text{contains but does not start at}  w\}$.
Then these have the following properties:

\begin{enumerate}
\item[1.]
the $A_w$'s are disjoint and $V=\cup_i A_i$,
\item[2.]
$B_w\subset \cup_{w'\in P(w)}A_{w'}$. In particular $B_r=\emptyset$.
\item[3.]
 if $i,j,j'$ are vertices of $T$ such that $i\in P(j), j'\in P(j), i\in P(j')$ and $x\in B_j\cap A_i$ 
then $x\in B_{j'}$,
\item[4.]
if $x\in A_i, y\in A_j$ and $x<y$ then $i\notin P(j)$,
\item[5.]
$e\in E$ if and only if $e\subset A_w\cup B_w$ for some $w$.
\end{enumerate}

This leads to the following definition of a {\it tree structure}.

\begin{definition}
\lbl{def.str}
Let $T=(W,F)$ be a tree, $V=\{x_1,\dots x_k\}$ be an ordered set and sets $A(w),B(w): w\in W$ 
satisfy the above properties 1.,2.,3.,4. Moreover let $|B_w|=b_w$.
 Then $(B_w:w\in W)$ is called a
$(T,V,(A_w, b_w:w\in W))-$ tree structure (tree structre for short). The set of all structures
is denoted by $\S(T,V,(A_w, b_w:w\in W))$.
\end{definition}

 What distinguishes tree structures are the sets $B_w$.  $B_w$ is an arbitrary 
subset of $A_{p(w)}\cup B_{p(w)}$ of $b_w$ elements. Hence we get the following observation.

\begin{proposition}
\lbl{prop.numstr}
The number of tree structures is $\prod_{r\neq w\in W}$ $\binom{a_{p(w)}+b_{p(w)}}{b_w}$.
\end{proposition}

\begin{remark}
\lbl{rem.s}
On the other hand, each tree structure on $T= (W,F,r), V=\{1,\dots, k\}$ determines a set
$T_v, 1\leq v\leq k$ of subtrees of $T$ so that $T_u\cap T_v\neq \emptyset$ if and only
if $u,v\in A_w\cup B_w$ for some $w\in W$, by reversing the construction of the tree
structure described above.
\end{remark}

\begin{definition}
\lbl{def.strdef}
Let $S$ be a tree structure, $v\in \{0,\dots, z-1\}^V$ and
$x\in A_w$ for some $w\in W$.
\begin{itemize}
\item
We denote by $G(S)$ the unique chordal graph with tree structure $S$ (see the remark above).
\item
We let $m(S,x)$ be the number of $y\in A_w\cup B_w$ such that $y<x$. Note that
$m(S,x)$ equals $b_w$ plus the number of elements of $A_w$ that are smaller than $x$.
Hence $m(S,x)$ does not depend on $S$ and we let $m(S,x)=m(x)$.
\item
We let $V(S,z)=\{v\in \{0,\dots, z-1\}^V;$
if $\{x,y\}\subset A_w\cup B_w$ for some $w$ then $v_x\neq v_y\}$. 
\item
We let $\defect(S,v,x)$ equal to the number of $y\in A_w\cup B_w$ such that 
$y< x$ and $v_y< v_x$.
\end{itemize}
\end{definition}

\begin{theorem}
\lbl{thm.str20}
$$
\sum_{S\in\S(T,V,(A_w, b_w:w\in W))}\prod_{w\in W}
\prod_{j=1}^{|A_w|}\frac{(b_w+j)_{q^{-1}}}{b_w+j} M_q(G(S),z)=
\prod_{r\neq w\in W}\binom{a_{p(w)}+b_{p(w)}}{b_w}
\prod_{x\in V} (z-m(x))_{q}.
$$
\end{theorem}

\subsection{A motivation from the quantum knot theory: colored Jones function}
\lbl{sub.mot}

 The motivation to the previous discussion comes from a study of the 
colored Jones function $J_n$ done jointly with Stavros Garoufalidis
in \cite{GL}. Colored Jones function is the quantum group invariant
of knots that corresponds to the $(n+1)$-dimensional irreducible representation
of $\mathfrak{sl}_2$. In \cite{GL} a new approach to  the colored Jones
function, based on the Bass-Ihara-Selberg zeta function of a graph, is presented.

Fix a generic planar projection $\K$ of an oriented knot with $r$ crossings. 
Let $c_i$ for $i=1,\dots,r$ denote an ordering of the crossings of $\K$. Then 
$\K$ consists of $r$ arcs $a_i$, which we label so that each arc $a_i$ ends 
at the crossing $i$. We will single out a specific arc of $\K$ which
we decorate by $\star$. Without loss of generality, we may assume
that the crossings of a knot appear in increasing order, when we walk
in the direction of the knot, and that the last arc is decorated by $\star$.

Given $\K$, we define a weighted directed 
graph $G_{\K}$ as follows:

\begin{definition}
\lbl{def.arcgraph}
The {\em arc-graph} $G_{\K}$ has $r$ vertices $1,\dots,r$, $r$ {\em blue
directed edges} $(v,v+1)$ ($v$ taken modulo $r$) and $r$ {\em red 
directed edges} $(u,v)$, where at the crossing $u$ the arc that crosses over 
is labeled by $a_v$.

The vertices of $G_{\K}$ are equipped with a sign, where $\sign(v)$ is the
sign of the corresponding crossing $v$ of $\K$, and the edges of $G_{\K}$ are
equipped with a weight $\b$, where the weight of the blue edge $(v,v+1)$ is 
$t^{-\sign(v)}$, and the weight of the red edge $(u,v)$ is $1-t^{-\sign(u)}$. 
Here $t$ is a variable.

Finally, $G_K$ denotes the digraph obtained by deleting vertex $r$
from $G_{\K}$.
\end{definition}

It is clear from the definition that from every vertex of $G_{\K}$, the blue
outdegree is $1$, the red outdegree is $1$, and the blue indegree is $1$.
It is also clear that $G_{\K}$ has a Hamiltonian cycle that consists of all
the blue edges. We denote by $e_i^b$ ($e_i^r$) the blue (red) edge {\em leaving
vertex} $i$.

\begin{example}
\lbl{ex.1}
For the figure 8 knot we have:
$$
\psdraw{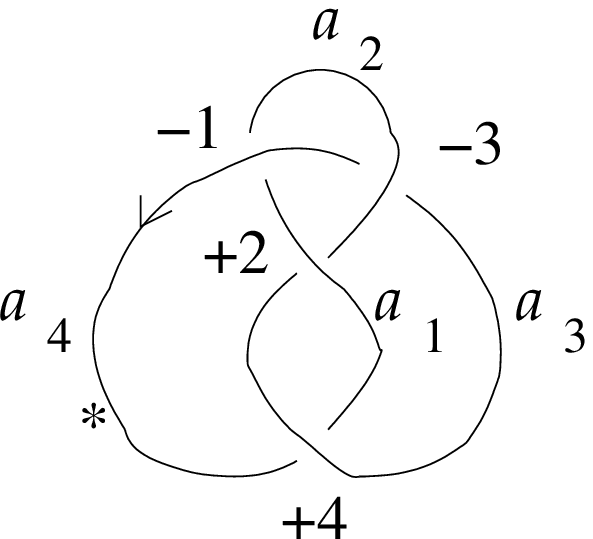}{1.5in}
$$
Its arc-graph $G_{\K}$ with the ordering and signs
of its vertices is given by 
$$
G_{\K}=\psdraw{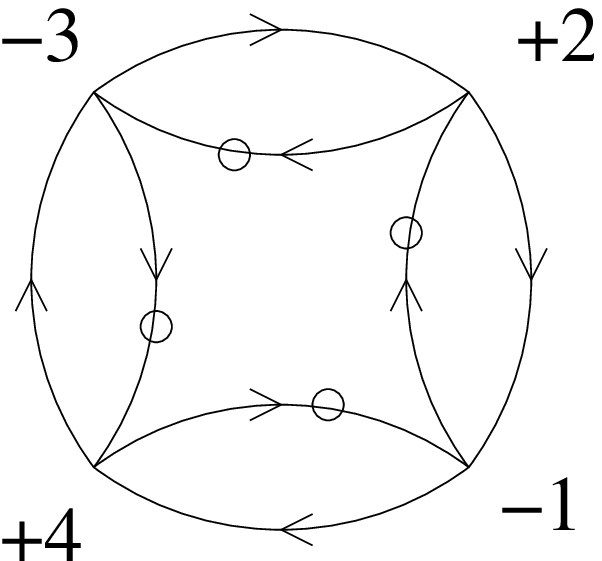}{1.3in} 
$$  
where the blue edges are the ones with circles on them.
\end{example}

\begin{definition}
\lbl{def.flow}  
A {\em flow} $f$ on a digraph $G$ is a function $f:
\mathrm{Edges}(G)\longto \BN$
of the edges of $G$ that satisfies the (Kirkhoff) {\em conservation law}
$$
\sum_{e \, \text{begins at} \, v} f(e)=\sum_{e \, \text{ends at} \, v}
f(e)
$$
at all vertices $v$ of $G$. Let $f(v)$ denote this quantity and let $\F(G)$ 
denote the {\em set of flows} of a digraph $G$.  

If $\beta$ is a weight function on
the set of edges of $G$ and $f$ is an flow on $G$, then
the {\em weight} $\b(f)$ of $f$ is given by
$\b(f)=\prod_e \b(e)^{f(e)}$, where $\b(e)$ is the weight of the edge $e$. 
\end{definition}

In order to express the Jones polynomial as a function of the reduced arc graph $G_K$, we 
need to add the following two structures, which may be read off
from the knot diagram.
\begin{itemize}
\item 
We associate in a standard way a {\em rotation} $\rot(e)$ to each edge $e$ of $G_K$; 
the exact definition is not relevant here; it may be found in \cite{GL}. 
\item
We linearly order the set of edges of $G_K$ terminating at vertex $v\in\{1,\dots, r-1\}$
as follows: if we travel along the arc of $K$ corresponding to vertex $v$, we 'see'
one by one the arcs corresponding to the starting vertices of red edges entering $v$:
this gives the linear order of the red edges entering $v$. Finally there is at most
one blue edge entering $v$, and we make it smaller than all the red edges entering $v$.
Let $P(e)$ denote the set of predecessors of an edge $e$ in the corresponding linear order.
\end{itemize}

With these decorations we define 
$$
\rot(f)=\sum_{e\in E} f(e)\rot(e),
\hspace{1cm}
\exc(f)=\sum_v \sign(v) f(e^b_v)(\sum_{e\in P(e^r_v)}f(e)),
\hspace{1cm} \d(f)=\exc(f)-\rot(f).  
$$

Let $\Sev(G)$ denote the set of all {\it admissible} subgraphs $C$ of $G$ such that each component
of $C$ is a directed cycle. Note that $\Sev(G)$ may be identified with a finite subset
of $\F(G)$ since the characteristic function of $C$ is a flow.

Let $\K$ be a knot projection. The {\em writhe} of $\K$, $\w(\K)$, is the sum of the signs of the
crossings of $\K$, and $\rot(K)$ is the {\em rotation number} of $K$, 
defined as follows: {\em smoothen} all crossings of $\K$, and consider the 
oriented circles that appear; one of them is special, marked by $\star$. The number of 
circles different from the special one whose orientation agrees with the 
special one, minus the number of circles whose orientation is opposite to the
special one is defined to be $\rot(K)$. 
We further let $\d(K,n)= 1/2 (n^2\w(\K)+n\rot(K))$, and $\d(K)= \d(K,1)$.

The following theorem appears in \cite{LW} (see also \cite{GL}).

\begin{theorem}
\lbl{thm.arcjones}
$$
J(\K)(t)= t^{\d(K)} \sum_{c \in 
\Sev(G_K)} t^{\d(c)} \b(c) .
$$
\end{theorem}

The colored Jones function equals to the Jones polynomial of 
a proper 'cabling' of the knot diagram. Using graph theory, this may be
described as follows.

\begin{definition}
\lbl{def.cabled2}
Fix a red-blue digraph $G$. Let $G^{(n)}$
denote the digraph with vertices $a^k_j$ for $k=1,\dots,r$ and $j=1,\dots,n$.
$G^{(n)}$ contains blue directed edges $(a^l_j,a^{l+1}_j)$ with weight
$t^{-\e n}$ (where $\e \in \{-1,+1\}$ is the sign of the crossing $l$)
for each
$l=1,\dots,r$ ($l+1$ considered modulo $r$) and $j=1,\dots,n$. Moreover,
if $(a_k,a_l)$ is a red directed edge of $G$, then $G^{(n)}$ contains
red edges $(a^k_i,a^l_j)$ for all $i,j=1,\dots n$ with weight
$t^{(j-1)}(1-t)$ resp. $t^{-(n-j)}(1-t^{-1})$, if the sign of the $i$ 
crossing is $-1$ resp. $+1$. 
\end{definition}

We will denote the set of admissible even subgraphs of $G^{(n)}$ by 
$\Sev_n(G)$. The following theorem appears in \cite{GL}.

\begin{theorem}
\lbl{thm.frst}
For every knot diagram $\K$ and every $n \in \BN$, we have
$$
J_n(\K)(t)=
t^{\d(K,n)} \sum_{c \in \Sev_n(G_K)} t^{\d(c)} \b(c). 
$$
\end{theorem}

Recall that for an integer $m$, we denote by 
$$
(m)_q=\frac{q^m-1}{q-1}
$$ 
the {\em quantum integer} $m$. This defines the 
{\em quantum factorial} and the {\em quantum binomial coefficients} by
$$
(m)_q!=(1)_q (2)_q \dots (m)_q
\qquad
\binom{m}{n}_q=\frac{(m)_q!}{(n)_q!(m-n)_q!}
$$
for natural numbers $m,n$ with $n \leq m$. We also define
$$
\mult_q(f)=\prod_v\binom{f(v)}{f(e_v^b)}_{q^{-\sign(v)}}.
$$

One of the key propositions of \cite{GL} is the following expression of the colored
Jones function (as a deformed zeta function of the reduced arc graph).

\begin{theorem}
\lbl{thm.main}
For oriented knot diagram $\K$ we have:
$$
J_n(\K) = t^{\d(K,n)}\sum_{f\in \F(G_K)}\mult_{t}(f)
t^{\d(f)} \prod_{v \in V_K}t^{-\sign(v)n f(e^b_v)}
\prod_{e \text{red}; t(e)=v} \prod_{j=0}^{f(e)-1} 
(1-t^{-\sign(s(e))(n-j-\sum_{e'<_v e} f(e))}).
$$
\end{theorem}

The proof of Theorem \ref{thm.main} is based on a rather complicated combinatorial construction. 
In this paper we present a curious interpretation of this construction as a defected q-chromatic
operator applied to a power series whose coefficients are linear combinations of chord diagrams.

\begin{definition}
\lbl{def.cd}
Given reduced knot diagram $K$ and a flow $f$ on $G_K$, we define:
\begin{itemize}
\item
A collection $I_1,\dots, I_p$ of intervals on $1,\dots, r-1$
is {\it relevant} (for $K,f$) if for each $1\leq v\leq r-1$, the number of intervals
starting at $v$ equals $\sum f(e); e$ red edge entering $v$, and the number of intervals
terminating at $v$ equals $f(e^r_v)$. 
We fix an order on the intervals of the relevant collection starting in the same vertex $v$,
according to $<_v$.
\item
Each relevant collection of intervals defines a set of {\it chord diagrams}. Chord diagram
means here set of chords of a line, with disjoint pairs of terminal vertices:
 
For each $1\leq v\leq r-1$ we introduce vertices 
$v_1,\dots, v_{o(v)}, v^1,\dots, v^{i(v)}$. We assume that all
the new vertices appear in the introduced order along a line.
For each interval we introduce a {\it chord} on this line. The chords corresponding to intervals
starting at $v$ will start at  $v_1,\dots, v_{o(v)}$, in agreement with the fixed ordering
of the intervals. The chords corresponding to intervals
terminating at $v$ will terminate at  $v^1,\dots, v^{i(v)}$, in {\it an arbitrary order}.
If $D$ is a resulting chord diagram, then we denote by $\deg(D)$ the number of 
chord diagrams obtained from the same relevant collection of intervals as $D$.
We assume that the chords in a chord diagram are ordered by their starting vertices.
\item
We denote by $\D(K,f)$ the set of chord diagrams obtained from a relevant collection
of intervals in this way.
\end{itemize}
\end{definition}

\begin{definition}
\lbl{def.intersection}
Let $D$ be a chord diagram. We define intersection graph of its chords $G(D)=(D,E(D))$ so that
the chords  of D form the set of vertices of $G(D)$, and two vertices form an edge if 
the corresponding chords intersect or one contains the other. 
\end{definition}

\begin{definition}
\lbl{def.chromdef}
Let $D$ be a chord diagram, $c\in D$ and $v\in V(G(D),n)$. We let 
\begin{itemize}
\item
We denote by $P(D,c)$ the set of chords $c'\in D$ which encircle the starting vertex
of $c$,
\item
$\defect_1(D,v,c)$ equals the number of chords $c'$ of $P(D,c)$ satisfying $v_{c'}< v_c$,
\item
We denote by $Q(D,c)$ the set of chords $c'\in D$ which encircle the terminal vertex $z$
of $c$ and at least one starting vertex of a chord after $z$,
\item
$\defect_2(D,v,c)$ equals the number of chords $c$ of $Q(D,c)$ satisfying  $v_{c'}< v_c$,
\item
$$
M^{\defect}_t(G(D),n)= 
\sum_{v\in V(G(D),n)}\prod_{c\in D} t^{v_c-\defect_1(D,v,c)-\defect_2(D,v,c)}.
$$
\end{itemize}
\end{definition}

\begin{remark}
\lbl{re.qdet}
We remark that $\defect_1$ is very close to the definition of sign for quantum determinants, 
first introduced by L. Fadeev, N. Reshetikhin and L. Takhtadjian in \cite{FRT}. Indeed,
Theorem \ref{thm.main} is used in \cite{GL} to give a non-commutative formula
for the colored Jones function.
\end{remark}

\begin{theorem}
\lbl{thm.ma2}
$$
J_n(K)(t)=
t^{\d(K,n)}\sum_{f\in \F(G_K)}Z_{n,f}(t)\sum_{D\in \D(K,f)}\deg(D)^{-1}M^{\defect}_t(G(C),n),
$$
where $Z_{n,f}(t)$ is a Laurent polynomial in $t$ with integer coefficients and parameters $f$ and $n$:
$$
Z_{n,f}(t)= 
t^{\d(f)}
t^{n(f^-_b-f^+_b)}(1-t)^{f^-_r}(1-t^{-1})^{f^+_r}
\prod_{e\in F^+_r} t^{-(n-1-|P(f,e)|)}
\prod_{v; \sign(v)= +} t^{f(e^r_v)f(e^b_v)}.
$$

Moreover for each $n$ only flows bounded by $n$ may contribute a non-zero to the RHS.

The terms $f^-_b, f^+_b, f^-_r, f^+_r, F^+_r$ are defined in \ref{sub.prbla} and the 
term $P(f,e)$ is defined in Definition \ref{def.prd}.
\end{theorem}

\subsection{Chord diagrams, Vassiliev invariants and Kontsevich integral}
\lbl{sec.Kon}
A describtion of the theory of Vassiliev knot invariants and weight systems
may be found in a seminal paper \cite{BN}.  It may be interesting
to explore a relationship of the formula of theorem \ref{thm.ma2} with
the Kontsevich integral expression for the colored Jones function.

\section{Proofs and Comments}
\lbl{sec.prf}

\subsection{The Principle of Inclusion and Exclusion And The Chromatic Polynomial}
\lbl{sub.chrom}

In 1932 Hassler Whitney \cite{Wh} deduced a formula for the chromatic polynomial of graphs
using the principle of inclusion and exclusion (PIE):

\medskip\noindent

If $A_1,...,A_n$ are finite sets, and if we let $\cap(A_i;i\in J)=A_J$ then 
$$|\cup(A_i;i=1,...,n)|=\sum_{k=1}^n (-1)^{k-1} \sum_{J\in {n\choose k}} |A_J|.$$

Let us present one of its folklore proofs, which uses binomial-type counting.

We use formula 

$$
\prod_{i=1}^n(1+x_i)= \sum_{I\subset \{1,\dots, n\}}\prod_{i\in I}x_i.
$$

Let $A=\cup_{1\leq i\leq n}A_i$ and let $f_i$ denote the characteristic function of $A_i$
in $A$. If $a\in A$ then $\prod_{i=1}^n (1-f_i(a))=0$, and so by the above formula
$$
\sum_{I\subset \{1,\dots, n\}}(-1)^{|I|}\prod_{i\in I}f_i(a)= 0.
$$

Summing these for each $a\in A$ we get
$$
0= \sum_{a\in A}\sum_{I\subset \{1,\dots, n\}}(-1)^{|I|}\prod_{i\in I}f_i(a)=
$$
$$
|A|+ \sum_{\emptyset\neq I\subset \{1,\dots, n\}}(-1)^{|I|}|\cap_{i\in I}A_i|
$$
since 
$$
\sum_{a\in A}\prod_{i\in I}f_i(a)= |\cap_{i\in I}A_i|.
$$

This is what we wanted to show.

\medskip\noindent

The chromatic polynomial of a graph $G=(V,E)$, denoted by $M(G,z)$, equals the number
of proper colorings of $G$ by $z$ or fewer colors. A {\it proper coloring}
is assigning one of the colors to each vertex of the graph in such a way that any two vertices 
which are joined by an edge are of different colors. 

Let $\{v_1,e_1,v_2,e_2,...,v_i,e_i,v_{i+1},...,e_n,v_{n+1}\}$ be a sequence
 such that each $v_j$ is a vertex of a graph $G$, each $e_j$ is an edge
 of $G$ and $e_j=v_jv_{j+1}$, and $v_i\neq v_j$ for $i<j$ except if
 $i=1$ and $j=n+1$. If also $v_1\neq v_{n+1}$ then $P$ is called {\it a path}
 of $G$. If $v_1=v_{n+1}$ then $P$ is called {\it a cycle} of $G$.
 In both cases the {\it length} of $P$ equals $n$. When no confusion arises we shall  
determine paths by listing their edges, namely  $P=(e_1,e_2,\dots,e_n)$. 
 A graph $G=(V,E)$ is {\it connected} if it has a path between any pair of vertices. 
If a graph is not connected then its maximum connected subgraphs are called
{\it connectivity components}. A subgraph of graph $G$ is {\it spanning} if its set of vertices
consists of all the vertices of $G$.

If $e\in E$ then let $A_e$ denote 
the set of the colorings with the property that the end-vertices are of the same color.
Then 
$$
M(G,z)= z^{|V|}-|\cup_ {e\in E} A_e|.
$$
If $G$ has $(p,s)$ (this is Birkhoff's symbol) spanning subgraphs of s edges
in p connectivity components, then by using PIE we get the well-known formula
for the chromatic polynomial:
$$
M(G,z) = \sum_{p,s}(p,s) (-1)^s z^p.
$$

Let $G=(V,E)$ be a graph. For $A \subset E$ let $r(A)=|V|-c(A)$,
 where $c(A)$ denotes the number of connectivity components of $G$. 
Then we can write 
$$
M(z)= z^{c(E)}(-1)^{r(E)}\sum_{A\subset E} (-z)^{r(E)-r(A)}(-1)^{|A|-r(A)}.
$$

This leads to Whitney rank generating function $R(G,u,v)$  defined by
$$
R(G,u,v)=\sum_{A\subset E} u^{r(E)-r(A)}v^{|A|-r(A)}.
$$
The Tutte polynomial has been defined by Tutte (\cite{T1}, \cite{T2}) as a 
minor modification of the Whitney rank generating function.
$$
T(G,x,y)=\sum_{A\subset E} (x-1)^{r(E)-r(A)}(y-1)^{|A|-r(A)}.
$$

In fact, both Whitney and Tutte polynomials are simply equivalent to a more straightforward
but less well-known generalization of the chromatic polynomial, the bichromatic polynomial
$$
B(G,a,b)= \sum_{A\subset E} a^{c(A)}b^{|A|}.
$$

\subsection{Geometric Summation and Quantum Binomial Formulas}
\lbl{sub.geom}
We all know the geometric summation formula
$$
\sum_{v_1,\dots, v_k=0}^{n-1}q^{\sum_i v_i} = [\frac{q^n-1}{q-1}]^k.
$$

The following quantum binomial formula leads to a well-known formula for the summation 
of  the products of distinct powers. We include a proof here in order
to keep the paper essentially self-contained. 

\begin{theorem}
\lbl{thm.qb}
$$
(a-z)(a-qz) \dots (a-q^{n-1}z)= \sum_{i=0}^n (-1)^i  {{n}\choose{i}}_q q^{i(i-1)/2}a^{n-i}z^i.
$$
\end{theorem}
\begin{proof}
We proceed by induction on $n$. It is easy to check the case $n=1$.
In the induction step assume the statement holds for $n$ and we want to prove it for $n+1$.
Let $y=qz$. We have
$$
(a-z)(a-qz) \dots (a-q^nz)= (a-z)(a-y) \dots (a-q^{n-1}y)=
(a-z)[\sum_{i=0}^n (-1)^i  {{n}\choose{i}}_q q^{i(i-1)/2}a^{n-i}y^i]=
$$
$$
\sum_{i=0}^n (-1)^i  {{n}\choose{i}}_q q^{i(i-1)/2}a^{n+1-i}z^i q^i+
\sum_{i=0}^n (-1)^{i+1} {{n}\choose{i}}_q q^{i(i-1)/2}a^{n-i}z^{i+1}q^i=
$$
$$
\sum_{i=0}^n (-1)^i  {{n}\choose{i}}_q q^{i(i-1)/2}a^{n+1-i}z^i q^i+
\sum_{i=1}^{n+1} (-1)^i {{n}\choose{i-1}}_q q^{i(i-1)/2}a^{n+1-i}z^i=
$$
$$
{{n}\choose{0}}_q a^{n+1} + (-1)^{n+1}{{n}\choose{n}}_q q^{n(n+1)/2}z^{n+1}+
$$
$$
\sum_{i=1}^n (-1)^i q^{i(i-1)/2}a^{n+1-i}z^i
[q^i{{n}\choose{i}}_q+{{n}\choose{i-1}}_q ]=
$$
$$
\sum_{i=0}^{n+1} (-1)^i  {{n+1}\choose{i}}_q q^{i(i-1)/2}a^{n+1-i}z^i
$$
since it may be observed directly that 
$$
q^i{{n}\choose{i}}_q+{{n}\choose{i-1}}_q = {{n+1}\choose{i}}_q.
$$
\end{proof}
Examining the coefficient of $z^k$ in the RHS, we get immediately

\begin{corollary}
\lbl{cor.dist}
$$
M_q(K_k,n)=k!{{n}\choose{k}}_q q^{k(k-1)/2}.
$$
\end{corollary}

\subsection{Proofs}
\lbl{sub.prd}

\begin{proof} (of theorem \ref{thm.qch})

If $A\subset E$ then let
$W(A,z)=\{v\in \{0,\dots, z-1\}^V$; if $\{i,j\}\in A$ then $v_i=v_j\}$.

The next considerations connect the PIE with the geometric series formula.
$$
M_q(G,z)=
\sum_{v\in \{0,\dots, z-1\}^V}q^{\sum_i v_i}- \sum_{v\in\cup_{e\in E}J_e}q^{\sum_i v_i},
$$
where $J_e, e=\{i,j\}\in E$, denotes the set of all vectors satisfying $v_i=v_j$.

By PIE this equals
$$
\sum_{A\subset E}(-1)^{|A|}\sum_{v\in\cap_{e\in A}J_e}t^{\sum_i v_i}=
\sum_{A\subset E}(-1)^{|A|}\sum_{v\in W(A,n)}t^{\sum_i v_i}=
$$
$$
\sum_{A\subset E}(-1)^{|A|}\prod_{W\in C(A)}\sum_{x\in\{0,\dots, z-1\}}q^{|W|x}=
\sum_{A\subset E} (-1)^{|A|} \prod_{W\in C(A)} (z)_{q^{|W|}}.
$$
\end{proof}

\begin{proof} (of theorem \ref{thm.qP})

We have
$$
P_q^k(G, J_e)= \sum_s q^{\sum_{v\in V}s(v)}\prod_{\{i,j\}\in E}
(1+ v_{ij}\delta(s(i), s(j)))=
$$
$$
\sum_s q^{\sum_{v\in V}s(v)}\sum_{A\subset E}\prod_{\{i,j\}\in A} v_{ij}\delta(s(i), s(j))=
$$
$$
\sum_{A\subset E}\sum_{s\in W(A,k)}q^{\sum_{v\in V}s(v)}
\prod_{\{i,j\}\in A} v_{ij}=
$$
$$
\sum_{A\subset E} \prod_{W\in C(A)} (k)_{q^{|W|}}\prod_{\{i,j\}\in A}v_{ij}.
$$
\end{proof}

\begin{proof} (of theorem \ref{thm.w})

Using the identity
 
$$
e^{x(s(i)s(j))}=cosh(x)+s(i)s(j)sinh(x),
$$
we have:
$$
\sum_{s}q^{\sum_{v\in V}s(v)}e^{\sum_{\{i,j\}\in E} J_{ij}s(i)s(j)}=
$$
$$
\sum_{s}q^{\sum_{v\in V}s(v)}\prod_{\{i,j\}\in E}[cosh( J_{ij})+s(i)s(j)sinh( J_{ij})]=
$$
$$
\prod_{\{i,j\}\in E}cosh( J_{ij})
\sum_{s}q^{\sum_{v\in V}s(v)}\prod_{\{i,j\}\in E}[1+s(i)s(j)th( J_{ij})]=
$$
$$
\prod_{\{i,j\}\in E}cosh( J_{ij})
\sum_{s}q^{\sum_{v\in V}s(v)}\sum_{A\subset E}\prod_{\{i,j\}\in A}
s(i)s(j)th( J_{ij})=
$$
$$
\prod_{\{i,j\}\in E}cosh( J_{ij})\sum_{A\subset E}\prod_{\{i,j\}\in A}th( J_{ij})
U(A),
$$
where 
$$
U(A)= \sum_{s}q^{\sum_{v\in V}s(v)}\prod_{\{i,j\}\in A}s(i)s(j).
$$
Theorem now follows from next Lemma \ref{lem.w}.
\end{proof}

\begin{lemma}
\lbl{lem.w}
Let $G= (V, E)$ be a graph. Then 
$$
\sum_{s}q^{\sum_{v\in V}s(v)}\prod_{\{i,j\}\in E}s(i)s(j)=
(q- q^{-1})^{o(E)}(q+ q^{-1})^{|V|- o(E)},
$$
where the first sum is over all functions $s$ from $V$ to $\{-1, 1\}$ and $o(E)$ denotes 
the number of vertices of $G$ of an odd degree.
\end{lemma}
\begin{proof}
First note that if $E$ is a cycle and $s$ arbitrary then $\prod_{\{i,j\}\in E}s(i)s(j)= 1$.
Hence, we can delete from $G$ any cycle without changing the LHS 
$$
\sum_{s}q^{\sum_{v\in V}s(v)}\prod_{\{i,j\}\in E}s(i)s(j).
$$
This reduces the proof to the case that $E$ is acyclic. If $E$ is a path, then it follows
from the observation above that for $s$ arbitrary, $\prod_{\{i,j\}\in E}s(i)s(j)= 1$
if and only if $s$ is constant on the end-vertices of $E$. Hence, we can delete from $E$
any maximal path and replace it by the edge between its end-vertices, without changing
the LHS. Hence it suffices to prove the proposition for the case that each component
of $G$ contains at most one edge. This is however simply true.

\end{proof}

\subsection{Proof of theorem \ref{thm.str20}}
\lbl{sub.thm2}

We first deduce a formula for a modified q-chromatic function.
Recall the definition of a tree structure $S$ for a chordal graph $G$,
and note that $V(S,z)= V(G,z)$.

\begin{proposition}
\lbl{prop.str2}
Let $G$ be a chordal graph and $S$ its tree structure. Then
$$
\sum_{v=(v_1\dots v_k)\in V(G,z)}q^{\sum_i v_i-\defect(S,v,i)}= 
\prod_{i=1}^k (z-m(i))_{q}.
$$
\end{proposition}
\begin{proof}

The basis for the calculation is the following Claim.

{\bf Claim.}
Fix numbers $v_1,\dots, v_{k-1}$ between $0$ and $z-1$ so that no edge of $G$
receives two equal numbers.
Then 
\begin{itemize}
\item
$$
\sum_{v_k:v=(v_1,\dots, v_k)\in V(G,z)} 
t^{v_k-\defect(v,k)}= A-B+C,
$$
where
$A=\sum_{v_k:v\in V(G,z)} t^{v_k}$, $B=\sum_{i=1}^{m(k)} t^{z-i}$, and
$C=\sum_{\{i,k\}\in E(G)} t^{v_i}$.
\item
$A+C= \sum_{0\leq j\leq z-1}t^z$ and $A-B+C= \frac{1-t^{z-m(k)}}{1-t}$.
\end{itemize}

{\bf Proof of Claim.}
Note that the second part simply follows from the first one. 

Let $v'_1 < \dots < v'_{m(k)}$ be a reordering of $\{v_i;\{i,k\}\in E(G)\}$.
We may write $v'_1=z-i_1, \dots,v'_{m(k)}=z-i_{m(k)}, 1\leq i_{m(k)}<\dots <i_1$. 
The LHS becomes

$$
(\sum_{v_k:v\in V(G,z)} t^{v_k})-t^{z-i_1+1}-\dots -t^{z-i_2-1}-t^{z-i_2+1}-\dots -
t^{z-i_{m(k)}-1}-t^{z-i_{m(k)}+1}- \dots -t^{z-1}+ t^{z-i_1}+\dots +
$$
$$
t^{z-i_2-2}+t^{z-i_2-1}+\dots +t^{z-m(k)-1}].
$$
This equals to the RHS of the equality we wanted to show. The Proposition simply follows
from the Claim.

\end{proof}

The proof of Proposition \ref{prop.str2} yields the following

\begin{proposition}
\lbl{prop.str.10}
Let $S$ be a structure. Then
$$
\sum_{v\in V(S,z)}\prod_{x\in V} q^{v_x-\defect(S,v,x)}= \prod_{x\in V} (z-m(x))_{q}.
$$
\end{proposition}

Hence 
$$
\sum_{v\in V(S,z)}\prod_{x\in V} q^{v_x-\defect(S,v,x)}
$$
is invariant for $S\in\S(T,V,(A_w, b_w:w\in W))$. Note that the same is not true for the non-defected version: path of three edges and star of three
edges, with their tree being the path, provide a contraexample.

\begin{proof} (of Theorem \ref{thm.str20})

$$
\sum_{S\in\S(T,V,(A_w, b_w:w\in W))}\sum_{v\in V(S,z)}\prod_{x\in V} q^{v_x-\defect(S,v,x)}=
\sum_{S=(B_w:w\in W)}\sum_{v\in V(S,z)}\prod_{x\in V} q^{v_x-\defect(S,v,x)}=
$$
$$
\sum_{v_x; x\in A_r}\sum_{B_w; p(w)=r}\prod_{x\in A_r}q^{v_x-\defect(S,v,x)}
[\sum_{(B_w: p(w)\neq r)}\sum_{v_x; x\notin A_r, v\in V(S,z)}
\prod_{x\in V-A_r} q^{v_x-\defect(S,v,x)}]=
$$
$$
\sum_{v(r)'}\sum_{B_w; p(w)=r}(|A_r|)!_{q^{-1}}\prod_{x\in A_r}q^{v(r)'_x}
[\sum_{(B_w: p(w)\neq r)}\sum_{v_x; x\notin A_r, v\in V(S,z)}
\prod_{x\in V-A_r} q^{v_x-\defect(S,v,x)}],
$$
where the second sum is over all vectors $v(r)'=(v(r)'_x; x\in A_r)$ so that
$v(r)'_x > v(r)_y$ for $x<y$. In the above equality we used
$$
\sum_{\pi}q^{\text{number}(i,j); i<j, \pi(i)< \pi(j)} = (n)!_q.
$$
This further equals
$$
\sum_{B_w; p(w)=r}(|A_r|)!_{q^{-1}}\sum_{v(r)'}\prod_{x\in A_r}q^{v(r)'_x}
[\sum_{(B_w: p(w)\neq r)}\sum_{v_x; x\notin A_r, v\in V(S,z)}
\prod_{x\in V-A_r} q^{v_x-\defect(S,v,x)}]=
$$
$$
\sum_{B_w; p(w)=r}\frac{(|A_r|)!_{q^{-1}}}{|A_r|!}\sum_{v_x; x\in A_r}\prod_{x\in A_r}q^{v_x}
[\sum_{(B_w: p(w)\neq r)}\sum_{v_x; x\notin A_r, v\in V(S,z)}
\prod_{x\in V-A_r} q^{v_x-\defect(S,v,x)}]=
$$
$$
\sum_{S\in\S(T,V,(A_w, b_w:w\in W))}\prod_{w\in W}
\prod_{j=1}^{|A_w|}\frac{(b_w+j)_{q^{-1}}}{b_w+j}
\sum_{v\in V(S,z)}\prod_{x\in V} q^{v_x}.
$$
Theorem now follows from Proposition \ref{prop.str.10}.
\end{proof}

\subsection{Categorification of flows: proof of Theorem \ref{thm.ma2}}
\lbl{sub.prbla}
Recall Theorem \ref{thm.frst}.
Each $c \in \Sev_n(G_K)$ projects to a flow on $G_K$.
An analysis of the contribution of each flow is obtained via
categorification of the flows and their multiplicities in \cite{GL}.
Next we briefly describe this.

Let $f$ be a flow on $G_K$. Let $F$ (resp. $F_r$) 
denote the multiset that contains each edge (resp. red edge) $e$ of 
$G_K$ with multiplicity $f(e)$. 

Let $F^+_r$ denote the set of all red edges of $F$ which leave a vertex with 
$+$ sign. Let $f^+_r=|F^+_r|$. Analogously we define  $F^-_r, \dots $.

If $e$ is an edge of $G_K$
then we let $F(e)\subset F$ be the set of all copies of $e$ in $F$, and fix
an {\em arbitrary} total order on each $F(e)$. We also denote by $t(e)$
the terminal vertex of $e$.

\begin{definition}
\lbl{def.conf}
Fix a flow $f$ on $G^\star_K$. A {\em flow configuration} of $f$
is a sequence 
$C=(C_1,\dots, C_{r-2})$ so that $C_1$ is a subset of $\{e\in F_r; e$ 
terminates in vertex $1\}$ of $f(e^b_1)$ elements and for each 
$2\leq i <r-1$,
$C_{i}$ is a subset of $C_{i-1}\cup \{e\in F_r; e$ terminates in vertex $i\}$ 
of $f(e^b_i)$ elements.
 \end{definition}

Let us denote by $\C(f)$ the set of all flow configurations of $f$. 

\begin{definition}
\lbl{def.AConf}
Let $f$ be a flow on $G_K$, $n>0$ a natural number, 
$C\in \C(f)$, and $v\in \{0,\dots, n-1\}^{F_r}$. 
We say that a pair $P= (C,v)$,  is {\em admissible} if,
for every two edges $e,e' \in F_r$ such that $v_e=v_{e'}$ and 
$e$ ends in vertex $i$ and $e'$ ends in vertex $j$ and $j \geq i$,
there exists an $l$, $i \leq l < j$ such that $e \not\in C_l$.
We denote the set of admissible flow configurations by $\AC(f)$.
\end{definition}

\begin{definition}
\lbl{def.prd}
Let $e'\in F(e)$. We define set $P(f,e')$ as follows:
if $e'_1\in F(e_1)$ then $e_1'\in P(f,e')$ if $e_1\in P(e)$ in $G_K$
or $e=e_1$ and $e'_1< e'$ in our fixed total order of $F(e)$.
\end{definition}

\begin{definition}
\lbl{def.okk}
Let $e \in F_r$. We define
\begin{itemize}
\item
$\defect_1(C,v,e)= |\{e'\in P(f,e): v_{e'}< v_{e} \}|$,
\item
$\defect_2(C,v,e)= |\{e'\in C_{d(e)}: v_{e'}< v_{e} \}|$,
where $d(e)$ is the biggest index such that $d(e)\geq t(e)$
and $e\not\in C_{d(e)}$.
\end{itemize}
\end{definition}

\begin{definition}
\lbl{def.Ifn}
If $v\in \{0,\dots,n-1\}^{F_r}$ then we define
$f^-_r(v)=\sum_{e\in F^-_r}v_e$ and we define $f^+_r(v)$ analogously. 
\end{definition}

The following theorem appears in \cite{GL}.

\begin{theorem}
\lbl{thm.catmm}
$$
J_n(K)(t)= t^{\d(K,n)}\sum_{f\in \F(G_K)}t^{\d(f)}
t^{n(f^-_b-f^+_b)}(1-t)^{f^-_r}(1-t^{-1})^{f^+_r}
$$
$$
 \prod_{e\in F^+_r} t^{-(n-1-|P(f,e)|)}
\prod_{i:\sign(i)=+}t^{f(e^r_v)f(e^b_v)}
\sum_{(C,v) \in \AC(f,n)}\prod_{e\in F_r} t^{v_e-\defect_1(C,v,e)-\defect_2(C,v,e)}.
$$
\end{theorem}

Finally we observe the relation of the flow structures and relevant collections of intervals.
Hence Theorem \ref{thm.catmm} implies Theorem \ref{thm.ma2}.

{\bf From flow structures to relevant collections of intervals.}
The formula of Theorem \ref{thm.catmm} may be interpreted in terms of chordal graphs.
The basic observation is that each flow structure is a $(T,V,(A_w, b_w:w\in W))-$ tree 
structure where $T$ is a path with vertices $1,\dots, r-1$ rooted at $1$, $B_{w+1}= C_w$
and $A_w= \{e\in F_r; e$ terminates in vertex $w\}$.
Now we recall that each chordal graph is the intersection graph of subtrees of a tree,
and this representation may be obtained from its tree structure (see subsection \ref{sub.chordal}).
However, if a tree structure of a graph  is indexed by a path, then it is the intersection graph of 
subpaths (intervals) of the path. This directly leads to the relevant collection of intervals, 
and to the proof of Theorem \ref{thm.ma2}.

\ifx\undefined\bysame
	\newcommand{\bysame}{\leavevmode\hbox
to3em{\hrulefill}\,}
\fi

\end{document}